\documentclass[a4paper]{article}
\usepackage{amsmath}
\usepackage{amssymb}
\usepackage{amsfonts}
\usepackage{theorem}
\usepackage{blkarray}
\usepackage{arydshln}
\usepackage{multirow}
\usepackage{graphicx}

\newcommand{\bigzerou}{%
\smash{\lower1.7ex\hbox{\bg 0}}}
\newtheorem{theorem}{Theorem}

\newtheorem{cor}{Corollary}

\newtheorem{lem}{Lemma}

\newcommand{\ba}{\begin{eqnarray}}
\newcommand{\ea}{\end{eqnarray}}
\newcommand{\ban}{\begin{eqnarray*}}
\newcommand{\ean}{\end{eqnarray*}}
\newcommand{\no}{\nonumber}

\newcommand{\e}{\epsilon}

\def\d{{\partial}}
\newcommand{\mapright}[1]{%
\smash{\mathop{%
\hbox to 1.0cm{\rightarrowfill}}\limits^{#1}}}
\newcommand{\mapleft}[1]{%
\smash{\mathop{%
\hbox to 1.3cm{\leftarrowfill}}\limits^{#1}}}
\allowdisplaybreaks[4]
\begin{document}
\title{
\begin{flushright}
  \begin{minipage}[b]{5em}
    \normalsize
    ${}$      \\
  \end{minipage}
\end{flushright}
{\bf Period Integrals (Givental's $I$-function) of Calabi-Yau Hypersurface in $CP^{N-1}$ and Intersection Numbers of Moduli Space of Quasimaps from $CP^{1}$ with Two Marked Points to $CP^{N-1}$ }}
\author{Masao Jinzenji${}^{(1)}$, Kohki Matsuzaka${}^{(2)}$ \\
\\ 
${}^{(1)}$ \it Department of Mathematics,  \\
\it Okayama University \\
\it  Okayama, 700-8530, Japan\\
\\
${}^{(2)}$\it Division of Mathematics, Graduate School of Science \\
\it Hokkaido University \\
\it  Kita-ku, Sapporo, 060-0810, Japan\\
\\
\it e-mail address: 
\it\hspace{0.0cm}${}^{(1)}$ pcj70e4e@okayama-u.ac.jp \\
\it\hspace{2.7cm}${}^{(2)}$ kohki@math.sci.hokudai.ac.jp } 
\maketitle
\begin{abstract}
In this paper, we derive the generalized hypergeometric functions (period integrals) used in mirror computation of Calabi-Yau hypersurface in $CP^{N-1}$ as generating functions of 
intersection numbers of the moduli space of quasimaps from $CP^{1}$ with two marked points to $CP^{N-1}$.   
\end{abstract}
\section{Introduction}
The aim of this paper is to clarify geometrical meaning of period integrals or generalized hypergeometric functions that are given as solutions of the following 
ordinary differential equation:
\ba
\biggl((\frac{d}{dx})^{N-1}-N\cdot e^{x}\cdot (N\frac{d}{dx}+(N-1))(N\frac{d}{dx}+(N-2))\cdots(N\frac{d}{dx}+1)\biggr)W(x)=0.
\label{eq2}
\ea
These were used in the mirror computation of genus $0$ Gromov-Witten invariants of Calabi-Yau hypersurface in $CP^{N-1}$ \cite{GMP,JN}.
By using Frobenius's method, the solutions are explicitly given as follows. 
\ba
W_{j}(x)&:=&\frac{\d^{j}}{\d \epsilon^{j}}\left(\sum_{d=0}^{\infty}\biggl(\frac{\prod_{r=1}^{Nd}(r+N\epsilon)}{\prod_{r=1}^{d}(r+\epsilon)^N}\biggr)e^{(d+\epsilon)x}\middle) \right|_{\epsilon=0}\no\\
         &=&\sum_{i=0}^{j}{j\choose i}x^{j-i}\biggl(\sum_{d=0}^{\infty}\frac{\d^{i}}{\d\epsilon^{i}}\left(\frac{\prod_{r=1}^{Nd}(r+N\epsilon)}{\prod_{r=1}^{d}(r+\epsilon)^N}\middle)\right|_{\epsilon=0}e^{dx}\biggr),\no\\
         &=:&\sum_{i=0}^{j}{j\choose i}x^{j-i}w_{i}(x) \;\;(j=0,1,\cdots,N-2).
\label{hgf}
\ea
On the other hand, geometrical study of mirror symmetry from the point of view of moduli space of quasimaps has been pursued by several groups of researchers \cite{CK,Jin1}.
Let us focus here the line of study in this direction pursued by our group. Here, we restrict our attention to the Calabi-Yau hypersurface in $CP^{N-1}$.
In \cite{Jin1}, we construct the moduli space $\widetilde{Mp}_{0,2}(N,d)$ of quasimaps (polynomial maps) 
of degree $d$ from $CP^{1}$ with two marked points to $CP^{N-1}$ and defined the intersection number $w({\cal O}_{h^{a}}{\cal O}_{h^{b}})_{0,d}$ ($a+b=N-3$ and $h$ is the hyperplane class in $H^{*}(CP^{N-1},{\bf C})$). This corresponds to B-model analogue of the genus $0$ and degree $d$ Gromov-Witten invariant $\langle{\cal O}_{h^{a}}{\cal O}_{h^{b}}\rangle_{0,d}$  of the Calabi-Yau hypersurface in $CP^{N-1}$. Let us introduce the following generating functions of these intersection numbers. 
\ba
&&w({\cal O}_{h^{a}}{\cal O}_{h^{b}})_{0}(x):=Nx+\sum_{d=1}^{\infty}w({\cal O}_{h^{a}}{\cal O}_{h^{b}})_{0,d}e^{dx},\no\\
&&\langle{\cal O}_{h^{a}}{\cal O}_{h^{b}}\rangle_{0}(t):=Nt+\sum_{d=1}^{\infty}\langle{\cal O}_{h^{a}}{\cal O}_{h^{b}}\rangle_{0,d}e^{dt}.
\ea
We have proved the following equalities:
\ba
&&w({\cal O}_{h^{N-3}}{\cal O}_{h^{0}})_{0}(x)=Nt(x)=N\frac{W_{1}(x)}{W_{0}(x)}=N(x+\frac{w_{1}(x)}{w_{0}(x)}),\no\\
&&\langle{\cal O}_{h^{a}}{\cal O}_{h^{b}}\rangle_{0}(t(x))=w({\cal O}_{h^{a}}{\cal O}_{h^{b}})_{0}(x).
\label{key2}
\ea
 The first equality tells us that $w({\cal O}_{h^{N-3}}{\cal O}_{h^{0}})_{0}(x)$ gives us the mirror map used in mirror computation of genus $0$ Gromov-Witten invariants of the 
 hypersurface. It was proved in \cite{Jin1}. The second equality tells that $w({\cal O}_{h^{a}}{\cal O}_{h^{b}})_{0}(x)$ is translated into generating function of the corresponding 
 Gromov-Witten invariants via the mirror map. Geometrical proof of it was given in \cite{Jin2}.  In \cite{S}, Saito gave explicit toric construction of $\widetilde{Mp}_{0,2}(N,d)$
 and showed that it is a compact orbifold. Moreover, he determined Chow ring of $\widetilde{Mp}_{0,2}(N,d)$. It is generated by $(d+1)$ generators $H_{0},H_{1},\cdots,H_{d}$ and relations 
 of the generators are given by, 
 \ba
(H_{0})^{N}=0,\;(H_{i})^{N}(2H_{i}-H_{i-1}-H_{i+1})=0\;\;(i=1,2,\cdots,d-1),\;(H_{d})^{N}=0. 
\ea
With these notations,  $w({\cal O}_{h^{a}}{\cal O}_{h^{b}})_{0,d}$ can be explicitly written as follows:
\ba
w({\cal O}_{h^{a}}{\cal O}_{h^{b}})_{0,d}=\int_{\widetilde{Mp}_{0,2}(N,d)}
(H_{0})^{a}\biggl(\frac{\prod_{j=1}^{d}e(H_{j-1}, H_{j})}{\prod_{j=1}^{d-1}(NH_{j})}\biggr)(H_{d})^{b},
\ea 
where $e(x,y)=\prod_{j=0}^{N}(jx+(N-j)y)$. This expression was already derived in \cite{Jin1} and was effectively used to prove the first equality of (\ref{key2}).  
Generalization of his explicit toric construction to the case of moduli space of quasimaps from $CP^{1}$ with two marked points to other toric manifolds have been 
given in \cite{JS2, Matsu}. 

Then we are naturally led to the question: \\
\\
`` {\bf Can we express the generalized hypergeometric function $W_{j}(x)$ in (\ref{hgf}) as generating function of intersection numbers 
of $\widetilde{Mp}_{0,2}(N,d)$?} ''\\
\\
Of course, we have already obtained some results on this question in the $N=5$ case, which were presented in Chapter 5 of \cite{Jin3}. Moreover, we 
can obtain many hints from Givental's theory of $I$-function and $J$-function \cite{giv}. In \cite{giv}, Givental suggested that the hypergeometric functions are closely related 
to the two point Gromov-Witten invariants:
\ba
\langle\sigma_{j}({\cal O}_{h^{N-3-j}}){\cal O}_{h^{0}}\rangle_{0,d}
:=\int_{\overline{M}_{0,2}(CP^{N-1},d)}(\psi_{1})^{j}\wedge \mbox{ev}_{1}^{*}(h^{N-3-j})\wedge c_{T}({\cal E}_{d}).
\label{givf}
\ea
 where $\overline{M}_{0,2}(CP^{N-1},d)$ is the moduli space of stable maps from genus $0$ stable curve with two marked points to $CP^{N-1}$,  
 $\psi_{1}$ is the Mumford-Morita-Miller class associated with the first marked point, $\mbox{ev}_{1}: \overline{M}_{0,2}(CP^{N-1},d)\rightarrow CP^{N-1}$ is the evaluation map at the first marked point and ${\cal E}_{d}$ 
 is the rank $Nd+1$ vector bundle that imposes condition that image of the stable curve lies inside the Calabi-Yau hypersurface. Moreover, we have known that 
 the element in Chow ring of  $\widetilde{Mp}_{0,2}(N,d)$ that corresponds to $\psi_{1}$ in the Chow ring of $\overline{M}_{0,2}(CP^{N-1},d)$ is given by $(H_{1}-H_{0})$ \cite{Kim, S}. 
 
 In this paper, we define the following intersection number on $\widetilde{Mp}_{0,2}(N,d)$:
 \ba
w(\sigma_{j}({\cal O}_{h^{a}}){\cal O}_{h^{b}})_{0,d}=\int_{\widetilde{Mp}_{0,2}(N,d)}
(H_{0})^{a}(H_{1}-H_{0})^{j}\biggl(\frac{\prod_{j=1}^{d}e(H_{j-1}, H_{j})}{\prod_{j=1}^{d-1}(NH_{j})}\biggr)(H_{d})^{b},
\label{wdefi}
\ea
 and prove the following theorem. 
\begin{theorem}
\ba
\frac{1}{N}w(\sigma_{j}({\cal O}_{h^{N-2-j}}){\cal O}_{h^{-1}}|{\cal O}_{h})_{0,2|1}&:=&\frac{d}{N}w(\sigma_{j}({\cal O}_{h^{N-2-j}}){\cal O}_{h^{-1}})_{0,2}+
\frac{1}{N}w(\sigma_{j-1}({\cal O}_{h^{N-1-j}}){\cal O}_{h^{-1}})_{0,2}\no\\
&=:&\frac{1}{j!}\frac{\d^{j}}{\d \epsilon^{j}}\left(\frac{\prod_{r=1}^{Nd}(r+N\epsilon)}{\prod_{r=1}^{d}(r+\epsilon)^N}\middle) \right|_{\epsilon=0}.
\label{mainth}
\ea
\label{main}
\end{theorem}
We have several remarks on the above theorem. Since we have insertion of ${\cal O}_{h^{-1}}$, we have to insert $\frac{1}{H_{d}}$ at the corresponding position in 
the formula (\ref{wdefi}). But this negative power cancels with $H_{d}$ in the polynomial $e^{N}(H_{d-1}, H_{d})=NH_{d-1}((N-1)H_{d-1}+H_{d})\cdots(H_{d-1}+(N-1)H_{d})NH_{d}$.
Hence the integrands used in evaluating $w(\sigma_{j}({\cal O}_{h^{N-2-j}}){\cal O}_{h^{-1}})_{0,2}$ and $w(\sigma_{j-1}({\cal O}_{h^{N-1-j}}){\cal O}_{h^{-1}})_{0,2}$ 
are polynomials in $H_{0},\cdots, H_{d}$\footnote{The operator insertion of ${\cal O}_{h^{-1}}$ was also considered by Zinger \cite{zinger} in defining the vector bundle ${\cal V}_{0}^{\prime}$ on $\overline{M}_{0,2}(CP^{N-1},d)$ in his notation.}. In (\ref{mainth}), we formally introduce the notation:
\ba
w(\sigma_{j}({\cal O}_{h^{N-2-j}}){\cal O}_{h^{-1}}|{\cal O}_{h})_{0,d}&:=&d\cdot w(\sigma_{j}({\cal O}_{h^{N-2-j}}){\cal O}_{h^{-1}})_{0,d}+w(\sigma_{j-1}({\cal O}_{h^{N-1-j}}){\cal O}_{h^{-1}})_{0,d}.
\label{hori}
\ea  
But we think that the intersection number $w(\sigma_{j}({\cal O}_{h^{N-2-j}}){\cal O}_{h^{-1}}|{\cal O}_{h})_{0,d}$ should be defined as three-pointed intersection number 
of the moduli space $\widetilde{Mp}_{0,2|1}(N,d)$ used in \cite{JS1}. Since we haven't explicitly determined Chow ring of $\widetilde{Mp}_{0,2|1}(N,d)$ and don't know well-defined 
Mumford-Morita-Miller class of the moduli space, the above notation is formal. If the notation $w(\sigma_{j}({\cal O}_{h^{N-2-j}}){\cal O}_{h^{-1}}|{\cal O}_{h})_{0,d}$ is justified as the corresponding intersection number of  $\widetilde{Mp}_{0,2|1}(N,d)$,  (\ref{hori}) corresponds to ``Hori's equation'' \cite{hori} for a three-pointed gravitational virtual structure constant:
\ba
w(\sigma_{j}({\cal O}_{h^{a}}){\cal O}_{h^{b}}|{\cal O}_{h})_{0,d}= d\cdot w(\sigma_{j}({\cal O}_{h^{a}}){\cal O}_{h^{b}})_{0,d}+w(\sigma_{j-1}({\cal O}_{h^{a+1}}){\cal O}_{h^{b}})_{0,d}.
\label{hori2}
\ea
Using this notation and Theorem \ref{main}, we can derive the following corollary by straightforward computation.
\begin{cor}
$W_{j}(x)$ is written as the following generating function of intersection numbers of  $\widetilde{Mp}_{0,2}(N,d)$.
\ban
\frac{1}{j!}W_{j}(x)&=&\sum_{d=0}^{\infty}\sum_{m=0}^{\infty}\sum_{k=0}^{\infty}\frac{x^m}{m!}\frac{1}{N}w(\sigma_{k}({\cal O}_{h^{N-2-j+m}}){\cal O}_{h^{-1}}|{\cal O}_{h})_{0,d}e^{dx}\\ 
                       &=:&\sum_{d=0}^{\infty}\frac{1}{N}w(\frac{e^{hx}h^{N-2-j}}{1-\tilde{\psi}_{1}},h^{-1}|h)_{0,d}e^{dx},
\ean
where $\tilde{\psi}_{1}$ denotes ``expected'' Mumford-Morita-Miller class of  $\widetilde{Mp}_{0,2|1}(N,d)$.
Then, we can also express Givental's $I$-function \cite{giv} as a generating function of intersection numbers of  $\widetilde{Mp}_{0,2}(N,d)$.
\ban
I(P,x):=\sum_{j=0}^{N-2}\frac{1}{j!}W_{j}(x)P^{j}&=&\sum_{j=0}^{N-2}\sum_{d=0}^{\infty}\frac{1}{N}w(\frac{e^{hx}h^{-j}P^{j}h^{N-2}}{1-\tilde{\psi}_{1}},h^{-1}|h)_{0,d}e^{dx}\\
&=&\frac{P^{N-2}}{N}\sum_{d=0}^{\infty}w(\frac{e^{hx}}{(1-\tilde{\psi}_{1})(1-\frac{h}{P})},h^{-1}|h)_{0,d}e^{dx}.
\ean
\end{cor}

\vspace{2em}
{\bf Acknowledgment} 
We would like to thank Prof. B. Kim and Prof. T. Ohmoto for valuable discussions. 
Our research is partially supported by JSPS grant No. 22K03289.  

\section{Proof of the Main Theorem }

Our proof of Theorem \ref{main} is based on the following formula derived in \cite{Jin1,S}:
\begin{align}
w(\sigma_{j}({\cal O}_{h^{a}}){\cal O}_{h^{b}})_{0,2} = &\frac{1}{(2 \pi \sqrt{-1})^{d+1}} \oint_{C_{0}} \frac{dz_{0}}{(z_{0})^{N}} \oint_{C_{1}} \frac{dz_{1}}{(z_{1})^{N}} \dots \oint_{C_{d}} \frac{dz_{d}}{(z_{d})^{N}} \no\\
& \times (z_{0})^{a} (z_{1} - z_{0})^{j} \prod_{l=1}^{d} e(z_{l-1},z_{l}) \prod_{l=1}^{d-1} \frac{1}{Nz_{l} (2z_{l} - z_{l-1} - z_{l+1})} \no\\
& \times (z_{d})^{b},
\label{residue}
\end{align}
where the operation $\frac{1}{2 \pi \sqrt{-1}} \oint_{C_{i}} dz_{i}$ means taking residues at $z_{i} = 0$ for $i = 0,d$ and at $z_{i} = 0, \frac{z_{i-1} + z_{i+1}}{2}$ for $i = 1 , \dots , d-1$. 
We introduce here a notation:
\[ w_{j,d} := w(\sigma_{j}({\cal O}_{h^{N-2-j}}){\cal O}_{h^{-1}})_{0,2} , \]
for brevity. 
$e(x,y)$ is given as, 
$$e(x,y)=\prod_{r=0}^{N}(rx+(N-r)y).$$
Note that $e(x,y) = e(y,x)$. For later use, we introduce,
\begin{align*}
a^{(l)}_i &:= \frac{1}{i!} \frac{\partial^{i}}{\partial \epsilon^i} \left( \prod_{r = 1}^{N} (N(l-1) + r + N \epsilon ) \middle) \right|_{\epsilon = 0} \\
&= \frac{1}{i!} \frac{\partial^{i}}{\partial \epsilon^i} \left( \prod_{r = N(l-1) + 1}^{Nl} (r + N \epsilon ) \middle) \right|_{\epsilon = 0} \quad (i=0 ,1,2, \dots ; l = 1 , \dots , d). 
\end{align*}
We can easily see, 
\begin{align*}
a^{(l)}_i := 0 \quad (i = N+1, N+2, \dots ;l = 1 , \dots , d).
\end{align*}
Before turning into the proof of Theorem \ref{main}, we prepare  several lemmas. 

\begin{lem}
\textit{For $l=1, \dots , d$, we can express $e(x,y)$ by using $a_{i}^{(l)}$ ($i = 1 , 2, \dots$) as follows.
\begin{equation}
e(x,y) = Nx \sum_{i = 1}^{\infty} a^{(l)}_{i} l^{i} \left( x - \frac{l-1}{l} y \right)^{i} (y-x)^{N-i}.
\label{l1}
\end{equation}
}
\end{lem}

\noindent
\textit{Proof}. \ Note that $\prod_{r = 1}^{N} (N(l-1) + r + N \epsilon )$ is expanded as, 
\[ \prod_{r = 1}^{N} (N(l-1) + r + N \epsilon ) = \sum_{i=0}^{\infty} a_i^{(l)} \epsilon^{i}. \]
By using this expression, 
we can derive (\ref{l1}) in the following way.
\begin{align*}
e(x,y) &= e(y,x) \\
        &= \prod_{r = 0}^{N} (ry + (N-r)x) \\
        &= Nx \prod_{r = 1}^{N} (r(y-x) + Nx) \\
        &= Nx(y-x)^{N} \prod_{r=1}^{N} \left( r + N \frac{x}{y-x} \right) \\
        &= Nx(y-x)^{N} \prod_{r=1}^{N} \left( N(l-1) + r + N \left( \frac{x}{y-x} -(l-1) \right) \right) \\
        &= Nx(y-x)^{N} \prod_{r=1}^{N} \left( N(l-1) + r + N \left( \frac{lx-(l-1)y}{y-x} \right) \right) \\
        &= Nx(y-x)^{N} \sum_{i=0}^{\infty} a_i^{(l)} \left( \frac{lx-(l-1)y}{y-x} \right)^{i} \\
        &= Nx \sum_{i=0}^{\infty} a_i^{(l)} (lx - (l-1)y)^{i} (y-x)^{N-i}. \\
        &= Nx \sum_{i=0}^{\infty} a_i^{(l)} l^{i} \left( x - \frac{l-1}{l} y \right)^{i} (y-x)^{N-i}. \ \square
\end{align*}

Next, we rewrite the r.h.s. of (\ref{mainth}).

\begin{lem}
\textit{
\begin{equation}
\frac{1}{j!} \frac{\d^{j}}{\d \epsilon^{j}}\left(\frac{\prod_{r=1}^{Nd}(r+N\epsilon)}{\prod_{r=1}^{d}(r+\epsilon)^N}\middle) \right|_{\epsilon=0} 
= \sum_{\substack{k_1 + \dots + k_d = j \\ k_1 , \dots , k_d \ge 0}} A_{k_1}^{(1)} \dotsb A_{k_d}^{(d)} ,
\label{l2}
\end{equation}
where
\[ A_{k}^{(l)} := \sum_{i = 0}^{k} a_{i}^{(l)} \binom{N+k-i-1}{k-i} \cdot \frac{(-1)^{k-i}}{l^{N+k-i}} . \]
}
\label{A}
\end{lem}

\noindent
\textit{Proof}. \ For $l = 1, \dots , d$, we set,
\[ F^{(l)} (\e ) := \prod_{r=N(l-1) +1}^{Nl} (r + N \e ) = \prod_{r = 1}^{N} (N(l-1) + r + N \epsilon ) = \sum_{i=0}^{\infty} a_i^{(l)} \epsilon^{i} , \] 
and 
\[ G^{(l)} (\e ) := (l + \e )^{-N}. \]
Then $(\d^{j} / \d \e^{j}) \{ F^{(1)} ( \e ) G^{(1)} ( \e ) \dotsb F^{(d)} ( \e ) G^{(d)} ( \e ) \} |_{\e = 0}$ is nothing but the l.h.s. of (\ref{l2}). 
By using the following identity:
\begin{align*}
G^{(l)} ( \e ) &= (l + \e )^{-N} \\ 
                &= \frac{1}{l^{N}} \left( 1 + \frac{\e}{l} \right) ^{-N} \\
                &= \frac{1}{l^{N}} \sum_{i = 0}^{\infty} \binom{-N}{i} \left( \frac{\e}{l} \right)^i \\
                &= \frac{1}{l^{N}} \sum_{i = 0}^{\infty} (-1)^i \binom{N+i-1}{N-1} \left( \frac{\e}{l} \right)^i \\
                &= \sum_{i = 0}^{\infty} \frac{1}{l^{N}}  \binom{N+i-1}{i} \left( - \frac{1}{l} \right)^i \e ^i ,
\end{align*}
we can rewrite $F^{(l)} ( \e ) G^{(l)} ( \e )$ as, 
\begin{align*}
F^{(l)} ( \e ) G^{(l)} ( \e ) &= \left( \sum_{k=0}^{\infty} a_k^{(l)} \epsilon^{k} \right) \left( \sum_{k = 0}^{\infty} \frac{1}{l^{N}} \binom{N+k-1}{i} \left( - \frac{1}{l} \right)^k \e ^k \right) \\
                                &= \sum_{k = 0}^{\infty} \left( \sum_{i=0}^{k} a_{i}^{(l)} \cdot \frac{1}{l^{N}} \binom{N+k-i-1}{k-i} \left( - \frac{1}{l} \right)^{k-i} \right) \e ^k \\
                                &= \sum_{k = 0}^{\infty} \left( \sum_{i = 0}^{k} a_{i}^{(l)} \binom{N+k-i-1}{k-i} \cdot \frac{(-1)^{k-i}}{l^{N+k-i}} \right) \e ^k \\
                                &= \sum_{k = 0}^{\infty} A_{k}^{(l)} \e ^k .
\end{align*}
Hence we obtain, 
\begin{align*}
F^{(1)} ( \e ) G^{(1)} ( \e ) \dotsb F^{(d)} ( \e ) G^{(d)} ( \e ) &= \left( \sum_{k_1 = 0}^{\infty} A_{k_1}^{(1)} \e ^{k_1} \right) \dotsb \left( \sum_{k_d = 0}^{\infty} A_{k_d}^{(d)} \e ^{k_d} \right) \\
&= \sum_{j=0}^{\infty} \left( \sum_{\substack{k_1 + \dots + k_d = j \\ k_1 , \dots , k_d \ge 0}} A_{k_1}^{(1)} \dotsb A_{k_d}^{(d)} \right) \e ^j. 
\end{align*}
This directly leads us to, 
\[ \frac{1}{j!} \frac{\d^{j}}{\d \epsilon^{j}}\left(\frac{\prod_{r=1}^{Nd}(r+N\epsilon)}{\prod_{r=1}^{d}(r+\epsilon)^N}\middle) \right|_{\epsilon=0} = \sum_{\substack{k_1 + \dots + k_d = j \\ k_1 , \dots , k_d \ge 0}} A_{k_1}^{(1)} \dotsb A_{k_d}^{(d)}. \ \square \]

The next two lemmas are useful in evaluating the residue integral in the r.h.s of (\ref{residue}),

\begin{lem}
\textit{Let $\alpha$ be a non-negative integer. For $l = 0 , \dots , d-2$, we have, 
\begin{align} 
&\frac{1}{2 \pi \sqrt{-1}} \oint_{z_{l} = \frac{l}{l+1} z_{l+1}} dz_{l} \frac{(z_{l+1} - z_{l})^{\alpha}}{\left( z_{l} - \frac{l}{l+1} z_{l+1} \right)^{\alpha + 1}} \cdot \frac{e(z_{l} , z_{l+1})}{Nz_{l} (2z_{l+1} - z_{l} - z_{l+2})} \no\\
&=(z_{l+1})^{N} \left( \frac{l+1}{l+2} \right)^{\alpha + 1} \sum_{s=0}^{\alpha} (l+2)^{s} A_{s}^{(l+1)} \frac{(z_{l+2} - z_{l+1})^{\alpha - s}}{\left( z_{l+1} - \frac{l+1}{l+2} z_{l+2} \right)^{\alpha - s + 1}}, 
\label{l3}
\end{align}
where the symbol $\frac{1}{2 \pi \sqrt{-1}} \oint_{z_{l} = \frac{l}{l+1} z_{l+1}} dz_{l}$ means taking a residue at $z_{l} = \frac{l}{l+1} z_{l+1}$. 
}
\end{lem}

\noindent
\textit{Proof.}
In this proof, we denote by $I$ the integral in the l.h.s of  (\ref{l3}). By using Lemma 1, the integrand of $I$ is rewritten as follows,
\begin{align*}
&\frac{(z_{l+1} - z_{l})^{\alpha}}{\left( z_{l} - \frac{l}{l+1} z_{l+1} \right)^{\alpha + 1}} \cdot \frac{e(z_{l} , z_{l+1})}{Nz_{l} (2z_{l+1} - z_{l} - z_{l+2})} \\
&= \frac{1}{\left( z_{l} - \frac{l}{l+1} z_{l+1} \right)^{\alpha + 1}} \cdot \frac{Nz_{l} \sum_{i=0}^{\infty} a_{i}^{(l+1)} (l+1)^{i} \left( z_{l} - \frac{l}{l+1} z_{l+1} \right)^{i} (z_{l+1} - z_{l})^{N + \alpha - i}}{Nz_{l} (2z_{l+1} - z_{l} - z_{l+2})} \\
&= \frac{1}{\left( z_{l} - \frac{l}{l+1} z_{l+1} \right)^{\alpha + 1}} \cdot \frac{\sum_{i=0}^{\infty} a_{i}^{(l+1)} (l+1)^{i} \left( z_{l} - \frac{l}{l+1} z_{l+1} \right)^{i} (z_{l+1} - z_{l})^{N + \alpha - i}}{2z_{l+1} - z_{l} - z_{l+2}}.
\end{align*}
Since we have, 
\begin{align*}
&\frac{\sum_{i=0}^{\infty} a_{i}^{(l+1)} (l+1)^{i} \left( z_{l} - \frac{l}{l+1} z_{l+1} \right)^{i} (z_{l+1} - z_{l})^{N + \alpha - i}}{2z_{l+1} - z_{l} - z_{l+2}} \Bigg|_{z_{l} = \frac{l}{l+1} z_{l+1}} \\
&= \frac{a_{0}^{(l+1)} \left( \frac{z_{l+1}}{l+1} \right)^{N + \alpha - i}}{\frac{l+2}{l+1} z_{l+1} - z_{l+2}} \\
&= \frac{\left( \frac{z_{l+1}}{l+1} \right)^{N + \alpha - i}}{\frac{l+2}{l+1} z_{l+1} - z_{l+2}} \cdot \prod_{r=1}^{N} (Nl + r) \neq 0,
\end{align*}
$I$ is explicitly evaluated as follows. 
\begin{align*}
I &= \frac{1}{\alpha !} \frac{\d ^{\alpha}}{\d z_{l}^{\alpha}} \Bigg|_{z_{l} = \frac{l}{l+1} z_{l+1}} \frac{\sum_{i=0}^{\infty} a_{i}^{(l+1)} (l+1)^{i} \left( z_{l} - \frac{l}{l+1} z_{l+1} \right)^{i} (z_{l+1} - z_{l})^{N + \alpha - i}}{2z_{l+1} - z_{l} - z_{l+2}} \\
&= \frac{1}{\alpha !} \sum_{p=0}^{\alpha} \binom{\alpha}{p} \frac{\d ^{p}}{\d z_{l}^{p}} \Bigg|_{z_{l} = \frac{l}{l+1} z_{l+1}} \left( \sum_{i=0}^{\infty} a_{i}^{(l+1)} (l+1)^{i} \left( z_{l} - \frac{l}{l+1} z_{l+1} \right)^{i} (z_{l+1} - z_{l})^{N + \alpha - i} \right) \\
& \quad \cdot \left( \frac{\d ^{\alpha - p}}{\d z_{l}^{\alpha - p}} \Bigg|_{z_{l} = \frac{l}{l+1} z_{l+1}} (2z_{l+1} - z_{l} - z_{l+2})^{-1} \right) \\
&= \frac{1}{\alpha !} \sum_{p=0}^{\alpha} \binom{\alpha}{p} \\
& \quad \cdot \sum_{i=0}^{\infty} a_{i}^{(l+1)} (l+1)^{i} \sum_{q=0}^{p} \binom{p}{q} \left( \frac{\d ^{q}}{\d z_{l}^{q}} \Bigg|_{z_{l} = \frac{l}{l+1} z_{l+1}} \left( z_{l} - \frac{l}{l+1} z_{l+1} \right)^{i} \right) \left( \frac{\d ^{p-q}}{\d z_{l}^{p-q}} \Bigg|_{z_{l} = \frac{l}{l+1} z_{l+1}} (z_{l+1} - z_{l})^{N + \alpha - i} \right) \\
& \quad \cdot \left( \frac{\d ^{\alpha - p}}{\d z_{l}^{\alpha - p}} \Bigg|_{z_{l} = \frac{l}{l+1} z_{l+1}} (2z_{l+1} - z_{l} - z_{l+2})^{-1} \right) . \\
\end{align*}
By basic calculus, we can easily see, 
\[ \frac{\d ^{q}}{\d z_{l}^{q}} \Bigg|_{z_{l} = \frac{l}{l+1} z_{l+1}} \left( z_{l} - \frac{l}{l+1} z_{l+1} \right)^{i} = \frac{i!}{(i-q)!} \delta_{0}^{i-q}, \]
\[ \frac{\d ^{p-q}}{\d z_{l}^{p-q}} \Bigg|_{z_{l} = \frac{l}{l+1} z_{l+1}} (z_{l+1} - z_{l})^{N + \alpha - i} = (-1)^{p-q} \frac{(N + \alpha - i)!}{(N + \alpha - i - p + q)!} \left( \frac{z_{l+1}}{l+1} \right)^{N + \alpha - i - p + q}, \]
\[ \frac{\d ^{\alpha - p}}{\d z_{l}^{\alpha - p}} \Bigg|_{z_{l} = \frac{l}{l+1} z_{l+1}} (2z_{l+1} - z_{l} - z_{l+2})^{-1} = (\alpha - p)! \left( \frac{l+2}{l+1} \right)^{- \alpha + p - 1} \left( z_{l+1} - \frac{l+1}{l+2} z_{l+2} \right)^{- \alpha + p - 1}, \]
where $\delta_{0}^{i-q}$ is Kronecker's delta. Therefore, we obtain,
\[ I = (z_{l+1})^{N} \sum_{p=0}^{\alpha} \sum_{i=0}^{p} a_{i}^{(l+1)} (l+1)^{-N+i + 1} (l+2)^{-\alpha + p - 1} (-1)^{p-i} \binom{N + \alpha - i}{p-i} \frac{(z_{l+1})^{\alpha - p}}{\left( z_{l+1} - \frac{l+1}{l+2} z_{l+2} \right) ^{\alpha - p + 1}}. \]
Now we expand $(z_{l+1})^{\alpha - p}$ in the following form:
\begin{align*}
(z_{l+1})^{\alpha - p} &= (l+2)^{\alpha - p} \left( \left( z_{l+1} - \frac{l+1}{l+2} z_{l+2} \right) + \frac{l+1}{l+2} (z_{l+2} - z_{l+1}) \right)^{\alpha - p} \\
&=(l+2)^{\alpha - p} \sum_{q=0}^{\alpha - p} \binom{\alpha - p}{\alpha - p -q} \left( z_{l+1} - \frac{l+1}{l+2} z_{l+2} \right)^{\alpha - p - q} \left( \frac{l+1}{l+2} \right)^{q} (z_{l+2} - z_{l+1})^{q}.
\end{align*}
Then $I$ is rewritten by, 
\begin{align*}
I &= (z_{l+1})^{N} \sum_{p=0}^{\alpha} \sum_{q=0}^{\alpha - p} \sum_{i=0}^{p} a_{i}^{(l+1)} (l+1)^{-N+i +q + 1} (l+2)^{-q- 1} (-1)^{p-i} \binom{N + \alpha - i}{p-i} \binom{\alpha - p}{\alpha - p -q} \\
& \quad \cdot \frac{(z_{l+2} - z_{l+1})^{q}}{\left( z_{l+1} - \frac{l+1}{l+2} z_{l+2} \right)^{q + 1}}.
\end{align*}
Note that the indices $(p,q,i)$ runs over the set $\{ (p,q,i) \in \mathbb{Z}^{3} \ | \ 0 \le p \le \alpha , 0 \le q \le \alpha - p, 0 \le i \le p \}$. 
Now we rename the indices as follows:
\[
\begin{cases}
s = q \\ 
t = i \\
u = p - i
\end{cases}
\]
Then $(s,t,u)$ runs over $\{ (s,t,u) \in \mathbb{Z}^{3} \ | \ 0 \le s \le \alpha , 0 \le t \le \alpha - s, 0 \le u \le \alpha - s - t \}$.  By using the identity,
\[ (-1)^{\alpha - s - t - u} \binom{\alpha - t - u}{\alpha - s - t - u} = \binom{-s-1}{\alpha - s - t - u}, \]
we can further rewrite $I$ into,  
\begin{align*}
I &= (z_{l+1})^{N} \sum_{s=0}^{\alpha} \sum_{t=0}^{\alpha - s} \sum_{u=0}^{\alpha - s - t} a_{t}^{(l+1)} (l+1)^{-N+s+t+1} (l+2)^{-s-1} (-1)^{\alpha - s - t} \binom{N + \alpha - t}{u} (-1)^{\alpha - s - t - u} \binom{\alpha - t - u}{\alpha - s - t - u} \\
&\quad \cdot \frac{(z_{l+2} - z_{l+1})^{s}}{\left( z_{l+1} - \frac{l+1}{l+2} z_{l+2} \right)^{s + 1}} \\
&= (z_{l+1})^{N} \sum_{s=0}^{\alpha} \sum_{t=0}^{\alpha - s} a_{t}^{(l+1)} (l+1)^{-N+s+t+1} (l+2)^{-s-1} (-1)^{\alpha - s - t} \\
&\quad \cdot \left( \sum_{u=0}^{\alpha - s - t} \binom{N + \alpha - t}{u} \binom{-s-1}{\alpha - s - t - u} \right) \\
&\quad \cdot \frac{(z_{l+2} - z_{l+1})^{s}}{\left( z_{l+1} - \frac{l+1}{l+2} z_{l+2} \right)^{s + 1}}.
\end{align*}
Then by using Chu--Vandermonde identity,
\[ \sum_{u=0}^{\alpha - s - t} \binom{N + \alpha - t}{u} \binom{-s-1}{\alpha - s - t - u} = \binom{N + \alpha - s - t - 1}{\alpha - s - t}, \]
we finally reach the r.h.s. of (\ref{l3}):
\begin{align*}
I &= (z_{l+1})^{N} \sum_{s=0}^{\alpha} \sum_{t=0}^{\alpha - s} a_{t}^{(l+1)} (l+1)^{-N+s+t+1} (l+2)^{-s-1} (-1)^{\alpha - s - t} \binom{N + \alpha - s - t - 1}{\alpha - s - t} \\
& \quad \cdot \frac{(z_{l+2} - z_{l+1})^{s}}{\left( z_{l+1} - \frac{l+1}{l+2} z_{l+2} \right)^{s + 1}} \\
&=(z_{l+1})^{N} \left( \frac{l+1}{l+2} \right)^{\alpha + 1} \sum_{s=0}^{\alpha} (l+2)^{\alpha -s} A_{\alpha - s}^{(l+1)} \frac{(z_{l+2} - z_{l+1})^{s}}{\left( z_{l+1} - \frac{l+1}{l+2} z_{l+2} \right)^{s + 1}} \\
&=(z_{l+1})^{N} \left( \frac{l+1}{l+2} \right)^{\alpha + 1} \sum_{s=0}^{\alpha} (l+2)^{s} A_{s}^{(l+1)} \frac{(z_{l+2} - z_{l+1})^{\alpha -s}}{\left( z_{l+1} - \frac{l+1}{l+2} z_{l+2} \right)^{\alpha -s + 1}}. \ \square
\end{align*}

\begin{lem}
\begin{align*}
&\frac{1}{2 \pi \sqrt{-1}} \oint_{z_{d-1} = \frac{d-1}{d} z_{d}} dz_{d-1} \frac{(z_{d} - z_{d-1})^{\alpha}}{\left( z_{d-1} - \frac{d-1}{d} z_{d} \right)^{\alpha + 1}} \cdot \frac{e(z_{d-1} , z_{d})}{Nz_{d-1}} \\
&= (z_{d})^{N} \cdot d^{\alpha} \cdot \sum_{l=0}^{\alpha} A_{l}^{(d)} \cdot (-d)^{- \alpha + l}.
\end{align*}
\end{lem}

\noindent
\textit{Proof.}
The proof goes in a similar manner as the one of  Lemma 3 except for use of the relation:
\[ \binom{N+a}{N} = \sum_{l=0}^{a} \binom{N+l-1}{N-1}. \]
 $\square$

\vspace{1cm}

\noindent
\textit{Proof of Theorem 1.}

\vspace{2mm}

The integrand in the r.h.s. of (\ref{residue}) is given by,
\begin{align*}
&\frac{1}{(z_{0})^{N} \dotsb (z_{d})^{N}} \cdot (z_{0})^{N-2-j} (z_{1} - z_{0})^{j} \cdot \frac{e(z_{0}, z_{1})}{Nz_{1} (2z_{1} - z_{0} - z_{2})} \dotsb \frac{e(z_{d-2},z_{d-1})}{Nz_{d-1} (2z_{d-1} - z_{d-2} -z_{d})} \cdot \frac{e(z_{d-1},z_{d})}{z_{d}} \\
&= \frac{1}{(z_{1})^{N} \dotsb (z_{d})^{N}} \cdot \frac{(z_{1} - z_{0})^{j}}{(z_{0})^{j+1}} \cdot \frac{e(z_{0}, z_{1})}{Nz_{0} (2z_{1} - z_{0} - z_{2})} \dotsb \frac{e(z_{d-2},z_{d-1})}{Nz_{d-2} (2z_{d-1} - z_{d-2} -z_{d})} \cdot \frac{e(z_{d-1},z_{d})}{Nz_{d-1}} \cdot \frac{N}{z_{d}}.
\end{align*}
By using Lemma 3, integration over $z_{0}$ results in the following integrand.
\begin{align*}
&\frac{1}{(z_{1})^{N} \dotsb (z_{d})^{N}} \cdot (z_{1})^{N} \left( \frac{1}{2} \right)^{j + 1} \sum_{s_{1} =0}^{j} 2^{s_{1}} A_{s_{1}}^{(1)} \frac{(z_{2} - z_{1})^{j -s_{1}}}{\left( z_{1} - \frac{1}{2} z_{2} \right)^{j -s_{1} + 1}} \\
&\cdot \frac{e(z_{1}, z_{2})}{Nz_{1} (2z_{2} - z_{1} - z_{3})} \dotsb \frac{e(z_{d-2},z_{d-1})}{Nz_{d-2} (2z_{d-1} - z_{d-2} -z_{d})} \cdot \frac{e(z_{d-1},z_{d})}{Nz_{d-1}} \cdot \frac{N}{z_{d}} \\
&= \frac{1}{(z_{2})^{N} \dotsb (z_{d})^{N}} \sum_{s_{1} =0}^{j} 2^{-j+s_{1} -1} A_{s_{1}}^{(1)} \frac{(z_{2} - z_{1})^{j -s_{1}}}{\left( z_{1} - \frac{1}{2} z_{2} \right)^{j -s_{1} + 1}} \cdot \frac{e(z_{1}, z_{2})}{Nz_{1} (2z_{2} - z_{1} - z_{3})} \\
&\cdot \frac{e(z_{2}, z_{3})}{Nz_{2} (2z_{3} - z_{2} - z_{4})} \dotsb \frac{e(z_{d-2},z_{d-1})}{Nz_{d-2} (2z_{d-1} - z_{d-2} -z_{d})} \cdot \frac{e(z_{d-1},z_{d})}{Nz_{d-1}} \cdot \frac{N}{z_{d}} \\
\end{align*}
Since this integrand is holomorphic at $z_{1} = 0$, we only have to take residue at $z_{1} = \frac{0 + z_{2}}{2} = \frac{z_{2}}{2}$. Thus integration over $z_{1}$ gives us,
\begin{align*}
&\frac{1}{(z_{2})^{N} \dotsb (z_{d})^{N}} \sum_{s_{1} =0}^{j} 2^{-j+s_{1} -1} A_{s_{1}}^{(1)} (z_{2})^{N} \left( \frac{2}{3} \right)^{j - s_{1} + 1} \sum_{s_{2} =0}^{j - s_{1}} 3^{s_{2}} A_{s_{2}}^{(2)} \frac{(z_{3} - z_{2})^{j - s_{1} -s_{2}}}{\left( z_{2} - \frac{2}{3} z_{3} \right)^{j  -s_{1} -s_{2} + 1}} \\
&\cdot \frac{e(z_{2}, z_{3})}{Nz_{2} (2z_{3} - z_{2} - z_{4})} \dotsb \frac{e(z_{d-2},z_{d-1})}{Nz_{d-2} (2z_{d-1} - z_{d-2} -z_{d})} \cdot \frac{e(z_{d-1},z_{d})}{Nz_{d-1}} \cdot \frac{N}{z_{d}} \\
&=\frac{1}{(z_{3})^{N} \dotsb (z_{d})^{N}} \sum_{s_{1} =0}^{j} A_{s_{1}}^{(1)} \sum_{s_{2} =0}^{j - s_{1}} 3^{-j + s_{1} + s_{2} -1} A_{s_{2}}^{(2)} \frac{(z_{3} - z_{2})^{j - s_{1} -s_{2}}}{\left( z_{2} - \frac{2}{3} z_{3} \right)^{j  -s_{1} -s_{2} + 1}} \cdot \frac{e(z_{2}, z_{3})}{Nz_{2} (2z_{3} - z_{2} - z_{4})} \\
&\cdot \frac{e(z_{3}, z_{4})}{Nz_{3} (2z_{4} - z_{3} - z_{5})} \dotsb \frac{e(z_{d-2},z_{d-1})}{Nz_{d-2} (2z_{d-1} - z_{d-2} -z_{d})} \cdot \frac{e(z_{d-1},z_{d})}{Nz_{d-1}} \cdot \frac{N}{z_{d}}.
\end{align*}
In the same manner, we integrate out $z_{2}, \dots , z_{d-2}$. Then we obtain,
\begin{align*}
w_{j,d} &= \frac{1}{(2 \pi \sqrt{-1})^{2}} \oint_{C_{d-1}} dz_{d-1} \oint_{C_{d}} dz_{d} \frac{1}{(z_{d})^{N}} \sum_{s_{1} =0}^{j} A_{s_{1}}^{(1)} \sum_{s_{2} =0}^{j - s_{1}} A_{s_{2}}^{(2)} \\ 
& \quad \dotsb \sum_{s_{d-1} = 0}^{j - s_{1} - \dotsb - s_{d-2}} d^{-j + s_{1} + \dotsb + s_{d-1} -1} A_{s_{d-1}}^{(d-1)} \frac{(z_{d} - z_{d-1})^{j - s_{1} \dotsb -s_{d-1}}}{\left( z_{d-1} - \frac{d-1}{d} z_{d} \right)^{j  -s_{1} \dotsb -s_{d-1} + 1}}  \cdot \frac{e(z_{d-1},z_{d})}{Nz_{d-1}} \cdot \frac{N}{z_{d}}.
\end{align*} 
Since this integrand is holomorphic at $z_{d-1} = 0$, we only have to take residue at $z_{d-1} = \frac{d-1}{d} z_{d}$. With the aid of Lemma 4, remaining integration goes as follows.
\begin{align*}
w_{j,d} &= \frac{1}{2 \pi \sqrt{-1}} \oint_{C_{d}} dz_{d} \sum_{s_{1} =0}^{j} A_{s_{1}}^{(1)} \sum_{s_{2} =0}^{j - s_{1}} A_{s_{2}}^{(2)} \\ 
& \quad \dotsb \sum_{s_{d-1} = 0}^{j - s_{1} - \dotsb - s_{d-2}} d^{-j + s_{1} + \dotsb + s_{d-1} -1} A_{s_{d-1}}^{(d-1)} \cdot d^{j - s_{1} - \dotsb - s_{d-1}} \cdot \sum_{s_{d} =0}^{j - s_{1} - \dotsb - s_{d-1}} A_{s_{d}}^{(d)} \cdot (-d)^{-j + s_{1} + \dotsb + s_{d}} \cdot \frac{N}{z_{d}} \\
&= \frac{N}{d} \sum_{s_{1} =0}^{j} A_{s_{1}}^{(1)} \sum_{s_{2} =0}^{j - s_{1}} A_{s_{2}}^{(2)} \\ 
& \quad \dotsb \sum_{s_{d-1} = 0}^{j - s_{1} - \dotsb - s_{d-2}} A_{s_{d-1}}^{(d-1)} \sum_{s_{d} =0}^{j - s_{1} - \dotsb - s_{d-1}} A_{s_{d}}^{(d)} \cdot (-d)^{-j + s_{1} + \dotsb + s_{d}} \\
&= \frac{N}{d} \left( - \frac{1}{d} \right)^{j} \sum_{r=0}^{j} \sum_{\substack{s_{1} + \dots + s_{d} = r \\ s_{1} , \dots , s_{d} \ge 0}} A_{s_{1}}^{(1)} \dotsb A_{s_{d}}^{(d)} (-d)^{r} \\
&= \frac{N}{d} \left( - \frac{1}{d} \right)^{j} \sum_{r=0}^{j} B_{r} (-d)^{r} ,
\end{align*}
where we used Lemma \ref{A} and introduced the notation: 
\[ B_{r} := \frac{1}{r!} \frac{\d^{r}}{\d \epsilon^{r}}\left(\frac{\prod_{l=1}^{Nd}(l+N\epsilon)}{\prod_{l=1}^{d}(l+\epsilon)^N}\middle) \right|_{\epsilon=0}.  \]
Therefore, the l.h.s. of (\ref{mainth}) reduces to, 
\begin{align*}
\frac{d}{N} w_{j,d} + \frac{1}{N} w_{j-1,d} &= \left( - \frac{1}{d} \right)^{j} \sum_{r=0}^{j} B_{r} (-d)^{r} + \frac{1}{d} \left( - \frac{1}{d} \right)^{j-1} \sum_{r=0}^{j-1} B_{r} (-d)^{r} \\
&= \left( - \frac{1}{d} \right)^{j} \sum_{r=0}^{j} B_{r} (-d)^{r} - \left( - \frac{1}{d} \right)^{j} \sum_{r=0}^{j-1} B_{r} (-d)^{r} \\
&= \left( - \frac{1}{d} \right)^{j} B_{j} (-d)^{j} \\
&= \frac{1}{j!} \frac{\d^{j}}{\d \epsilon^{j}}\left(\frac{\prod_{r=1}^{Nd}(r+N\epsilon)}{\prod_{r=1}^{d}(r+\epsilon)^N}\middle) \right|_{\epsilon=0},
\end{align*}
which completes the proof of Theorem \ref{main}. $\square$

\end{document}